\documentclass[german,english
              ,twoside
              ,10pt
              ]{amsart}
\usepackage{gt-sim}
\usepackage[latin1]{inputenc}
\usepackage{babel}
\usepackage{hyperref}
\usepackage{cas-math}
\usepackage{cas-paper}
\usepackage[arrow,cmtip,matrix]{xy}

\newcommand\w{w}
\newcommand\HHom{\Hom(C_{2r+1},K_{n+2})}
\DeclareMathOperator\bsd{bsd}
\DeclareMathOperator\esd{esd}
\newcommand\real{\abs}

\newenvironment{acknowledgements}{\subsection{Acknowledgements}}{}
\newcommand{\subclass}[1]{}

\begin{document}

\title[$\w_1^n(\Hom(C_{2r+1},K_{n+2}))=0$]{A short proof of 
$\w_1^n(\Hom(C_{2r+1},K_{n+2}))=0$ for all~$n$
\\and a graph colouring theorem 
by\\Babson and Kozlov}
\author{Carsten Schultz}
\address{Institut für Mathematik, MA 6-2\\
Technische Universität Berlin\\
D-10623 Berlin, Germany}
\email{carsten@codimi.de}
\thanks{This research was supported by the
\foreignlanguage{german}{Deutsche Forschungsgemeinschaft} within the
European graduate program ``Combinatorics, Geometry, and Computation''
(No.~GRK~588/2)}

\begin{abstract}
We show that the $n$-th power of the first Stiefel-Whitney class of the
\hbox{$\Z_2$-operation} on the graph complex $\HHom$ is zero, confirming a
conjecture by Babson and Kozlov.  This proves the strong form of their
graph colouring theorem, which they had only proven
for odd~$n$.  Our proof is also considerably simpler than their
proof of the weak form of the theorem, which is also
known as the Lovász conjecture.
\subclass{57M15; 05C15}
\end{abstract}

\maketitle

\section{Introduction}
A homomorphism from a graph~$G$ to a graph~$H$ is a function $f\colon
V(G)\to V(H)$ between their vertex sets 
that respects the edge relation, i.e.\ such that
$\set{f(u),f(v)}\in E(H)$ whenever $\set{u,v}\in E(G)$.  If the
vertices of~$H$ are identified with the vertices of the
simplex~$\Delta^{\abs{V(H)}-1}$, then functions $V(G)\to V(H)$
correspond to vertices of the cell complex~$\prod_{v\in
V(G)}\Delta^{\abs{V(H)}-1}$.  The subcomplex consisting of all cells
all of which vertices correspond to homomorphisms is denoted
by~$\Hom(G, H)$.  $\Hom$ is a functor, for its basic properties, as
well as background on the results mentioned in this introduction, we
refer the reader to \cite{babson-kozlov-i} and~\cite{kozlov-survey}.  
We will however try to be
as self-contained as possible in the proof of our main result,
\prettyref{thm:wv}.

The main result of~\cite{babson-kozlov-ii} is the following theorem.
$C_s$ denotes the circular graph with $s$~edges.
\begin{thm}\label{thm:bk}
Let $G$ be a graph, $r\ge1$ and $n\ge0$.  If $\Hom(C_{2r+1}, G)$ is
$(n-1)$-connected, then $G$ is not $(n+2)$-colourable.
\end{thm}
This had been conjectured by Lovász and is analogous to the following
theorem that represented his novel approach to graph colourings in his
proof of Kneser's conjecture~\cite{lovasz}.
\begin{thm}\label{thm:lovasz}
Let $G$ be a graph and $n\ge-1$.  If $\Hom(K_2, G)$ is
$n$-connected, then $G$ is not $(n+2)$-colourable.
\end{thm}
Lovász did originally not use the complex $\Hom(K_2,G)$, but another
graph complex, the neighbourhood complex, which is homotopy equivalent
to it.

The proofs of both theorems rely on the study of a free $\Z_2$-action
on the space $\Hom(C_{2r+1}, G)$ or $\Hom(K_2, G)$ respectively, which is
induced by the nontrivial automorphism of $K_2$ respectively an
automorphism of $C_{2r+1}$ flipping one edge.  We will always regard
these spaces as equipped with these actions.  Similarly, when a
$\Z_2$-action on a sphere is mentioned, this will always be given by
the antipodal map unless noted otherwise.

The basic idea of the proofs of both theorems is the same.  An
$(n+2)$-colouring of $G$ is just a graph homorphism $G\to K_{n+2}$ and
induces a map $\Hom(H, G)\to\Hom(H, K_{n+2})$, which, by functoriality
of~$\Hom$ with respect to its first argument, 
is a $\Z_2$-map for $H=K_2$ or $H=C_{2r+1}$.  If $\Hom(H,
G)$ is $(k-1)$-connected, then a $\Z_2$-map $\Sphere^k\to\Hom(H, G)$
and hence $\Sphere^k\to\Hom(H, K_{n+2})$ exists.  A study of the space
$\Hom(H, K_{n+2})$ then reveals conditions on $r$ and~$n$ under which
such a map does not exist.  In the case of \prettyref{thm:lovasz} this
amounts to the realisation that $\Hom(K_2,
K_{n+2})\homeo_{\Z_2}\Sphere^n$ (see \prettyref{rem:homk2kn}) 
and an invocation of the Borsuk-Ulam
Theorem.

Babson's and Kozlov's proof of \prettyref{thm:bk} uses different
arguments for odd and for even~$r$; both rely on elaborate cohomology
calculations for the cell complex $\HHom$.  For even~$n$ they show that a
$\Z_2$-map $\HHom\to\Sphere^n$ exists such that the image of a
generator of $H^n(\Sphere^n)$ in $H^n(\HHom)$ under the induced map is
of order~$2$.  Hence the degree of the composition of this map with a
map $\Sphere^n\to\HHom$ is zero.  A very elegant and simple proof of
this fact, avoiding all cohomology calculations, has recently been
found by \u Zivaljevi\'c~\cite{z-groupoids}.  A $\Z_2$-map
$\Sphere^n\to\Sphere^n$, however, has odd degree.

For odd~$n$, Babson and Kozlov use the same $\Z_2$-map and the map
between quotient spaces
$\HHom/\Z_2\to\R P^n$ it induces and show that the image of the
generator of~$H^n(\R P^n;\Z_2)$ in $H^n(\HHom/\Z_2;\Z_2)$ is zero.  Thus
for odd~$n$ they proved the following:
\begin{thm}\label{thm:wv}
Let $r\ge 1$, $n\ge0$.  Then $\w_1^n(\HHom)=0$.
\end{thm}
\begin{sloppypar}
Here $w_1$ denotes the first Stiefel-Whitney class of the free $\Z_2$-space
$\HHom$ (which can be regarded as a principal $\Z_2$-bundle over
$\HHom/\Z_2$, which determines a real line bundle since $\Z_2\isom O(1)$).
The statement $\w_1^0(X)=0$ is to be read as $H^0(X/\Z_2;\Z_2)=0$.
\end{sloppypar}

Babson and Kozlov have conjectured that \prettyref{thm:wv} holds for
all~$n$ \cite[Conj.~2.5]{babson-kozlov-ii}, and the objective of this
article is to prove this, by an argument that is valid for all~$n$.

From \prettyref{thm:wv} using the naturality of $\w_1$ with respect to
$\Z_2$-maps and $\w_1^n(\Sphere^n)\ne0$ one immediately deduces
\prettyref{thm:bk}.  Using only the naturality of $w_1$
we can also state the following stronger version which has been announced
by Babson and Kozlov in \cite{babson-kozlov-era}.
\begin{thm}
Let $G$ be a graph, $r\ge1$ and $n\ge0$.  If $\w_1^n(\Hom(C_{2r+1},G))\ne0$, 
then $G$ is not $(n+2)$-colourable.
\end{thm}
This version is not only stronger, it also stresses the importance of
the $\Z_2$-action and may be better suited to actual computations.

\subsection{Outline of the proof of \prettyref{thm:wv}}
We will construct an equivariant map of 
free $\Z_2$-spaces $f\colon\Hom(C_{2r+1}, K_{n+2})\to X_{n,r}$
where $X_{n,r}$ is a product of $r$ spheres of dimension $n$ with a
certain operation of~$\Z_2$. The map $f$ misses a certain $\Z_2$-invariant
subspace $A_{n,r}\subset X_{n,r}$.  An elementary computation shows
that the $n$-th power of the Stiefel-Whitney class~$\w_1$ of
$X_{n,r}\wo A_{n,r}$, and hence of $\Hom(C_{2r+1}, K_{n+2})$, is zero.

\subsection{Subsequent developments}
The proof presented here was the first proof of \prettyref{thm:wv} and
the first simple proof of \prettyref{thm:bk}.  Later on,
in~\cite{hom-loop} the more general result that $\w_1^n(\Hom(C_{2r+1},
G))=0$ holds for all graphs~$G$ with $\w_1^{n+1}(\Hom(K_2, G))=0$ has been
proven.  In~\cite{dmitry-cobounding} a proof of \prettyref{thm:wv} is
given that proceeds similarly to the proof in this work but uses
calculations at the cochain level instead of our calculations in the
homology groups of products of spheres.

\begin{acknowledgements}
The construction presented here has been inspired by a conjecture of
Péter Csorba \cite[Conj.~4.8]{csorba-thesis} 
that I learned about from Frank Lutz and that has lead to their
work~\cite{csorba-lutz} on $\Hom$-complexes that are manifolds.  That
conjecture states that $\Hom(C_5,K_{n+2})$ is homeomorphic to the
Stiefel manifold $\set{(x,y)\in\Sphere^n\times\Sphere^n\colon \langle
x,y\rangle=0}$ and has been proved in~\cite{c5}.
\end{acknowledgements}

\section{The map $\HHom\to X_{n,r}\wo A_{n,r}$}
\label{sec:X}
Let $r\ge1$, $n\ge0$.  

\begin{nota}
We denote the vertices of $C_{2r+1}$ by $v_0,\dots,v_{2r}$, in some
cyclic order.  The barycentric subdivision of the complex $\HHom$ can
be described as the order complex of the poset
\begin{multline*}
\big\{\phi\colon V(C_{2r+1})\to \power(V(K_{n+2}))\wo\set\emptyset
    \colon\\ 
    \text{$\phi(v_i)\intersect\phi(v_{i+1})=\emptyset$, $0\le i<2r$,
          $\phi(v_0)\intersect\phi(v_{2r})=\emptyset$}\big\}
\end{multline*}
ordered by 
$\phi\le\phi'\iff\text{$\phi(v_i)\subset\phi'(v_i)$ for all~$i$}$.
A vertex of this subdivision can be regarded as a multi-colouring
of~$C_{2r+1}$, in which each vertex is assigned a non-empty set of
colours, i.e.\ vertices of~$K_{n+2}$, and adjacent vertices are
assigned disjoint sets.\sloppy
We will always work with this subdivision.
The $\Z_2$-operation on $\HHom$ induced by the automorphism of
$C_{2r+1}$ that flips the edge $\set{v_0,v_{2r}}$ can be described
explicitly by $(\tau \phi)(v_i)=\phi(v_{2r-i})$, 
where $\tau$ denotes the generator of $\Z_2$.  
\end{nota}

We will compare $\HHom$ to the following space.

\begin{defn}
We define the space $X_{r,n}$ by
\begin{equation*}
X_{r,n}\deq \prod_{i=0}^{2r-1}\Sphere^n=(\Sphere^n)^{2r}
\end{equation*}
and a $\Z_2$-operation on it by
\begin{align*}
(\tau x)_0&\deq -x_0,\\
(\tau x)_i&\deq x_{2r-i},\qquad\text{$0<i<2r$.}
\end{align*}
\end{defn}
We will define an
equivariant map $f=(f_i)_i\colon \Hom(C_{2r+1}, K_{n+2})\to X_{n,r}$
as follows.  The inclusion of $K_2$ to the edge $\set{v_0,v_{2r}}$ of
$C_{2r+1}$ induces a map \[\HHom\to\Hom(K_2,K_{n+2})\homeo\Sphere^n.\]
This will be $f_0$. The inclusion of the vertex $v_{i}$ induces a map
\[\HHom\to\Hom(K_1,K_{n+2})\homeo\Disk^{n+1}.\]  Since $v_{2i}$ is not
isolated, the image of this map is contained in (actually equal to)
$\dd\Disk^{n+1}=\Sphere^n$.  This will be $f_i$ for $i>0$.  We now carry out
this construction in more detail.

We identify the vertices of $K_{n+2}$ with the vertices of the
$(n+1)$-simplex $\Delta^{n+1}$.  The barycentric subdivision of the
boundary of this simplex, $\bsd(\dd \Delta^{n+1})$, is the order
complex of the set $\set{A\colon A\subset V(K_{n+2}), A\ne\emptyset,
\compl A\ne\emptyset}$ ordered by inclusion.  The map sending a subset
of $V(K_{n+2})$ to its complement is a simplicial map of this complex
to itself, and we choose a triangulation $t\colon
\real{\bsd(\dd\Delta^{n+1})}\to\Sphere^n$ of the $n$-sphere that transports
this map to the antipodal map of the sphere.

\begin{defn}
For $0<i<2r$ we define $f_i\colon\HHom\to\Sphere^n$ as the composition
of the simplical map 
\begin{align*}
\HHom&\to\bsd(\dd\Delta^{n+1}),\\
\phi&\mapsto\phi(v_{i})
\end{align*}
with the triangulation~$t$.
\end{defn}

\begin{lem}\label{lem:Ai}
For $0<i<2r-1$ and any point $y\in\HHom$ we have $f_i(y)\ne f_{i+1}(y)$.
\end{lem}

\begin{proof}
Let $y\in\HHom$.  The point~$y$ is contained in a simplex defined by a
chain $\phi_0<\dots<\phi_s$.  The points $f_i(y)$ and $=f_{i+1}(y)$
are cointained in the simplices defined by the chains
$\phi_0(v_i)\subset\dots\subset\phi_s(v_i)$ and
$\phi_0(v_{i+1})\subset\dots\subset\phi_s(v_{i+1})$ respectively.  Since
$\phi_s(v_i)\intersect\phi_s(v_{i+1})=\emptyset$, these simplices have
no vertex in common and hence an empty intersection.
\end{proof}

For a poset $P$ we can consider its corresponding
\emph{interval poset} 
$\set{(p,q)\colon p\le q}$
ordered by $(p,q)\le(p',q')\iff(p\le p')\land(q\ge q')$, the order
complex of which is the \emph{edge subdivision}
$\esd(\Delta(P))$ of the order complex~$\Delta(P)$ of~$P$.  It has as
vertices the barycentres of $0$-simplices and $1$-simplices
of~$\Delta(P)$, see~\cite{walker-posets}.  
To see that this is actually a subdivision, note that
the edge subdivision of the simplex spanned by $p_0<\dots<p_s$ can
recursively be defined as the join of the barycentre of $\simp{p_0,p_s}$
with the union of the edge subdivisions of $\simp{p_0,\dots,p_{s-1}}$
and $\simp{p_1,\dots,p_s}$, which agree on $\simp{p_1,\dots,p_{s-1}}$.

The triangulation~$t$ induces a triangulation
$t'\colon\esd(\bsd(\dd\Delta^{n+1}))\to\Sphere^n$.

\begin{rem}\label{rem:homk2kn}\sloppy
The barycentric subdivision of $\Hom(K_2,K_{n+2})$ is the order
complex of $\set{(A,B)\colon A,B\subset V(K_{n+2}), A\intersect
B=\emptyset, A,B\ne\emptyset}$.  Via $(A,B)\mapsto(A,\compl B)$ this
is isomorphic to $\esd(\bsd(\dd\Delta^{n+1}))$.  Using $t'$ we have a
homeomorphism $\Hom(K_2,K_{n+2})\homeo\Sphere^n$.  This homeomorphism
transports the automorphism of $\Hom(K_2,K_{n+2})$ given by
$(A,B)\mapsto(B,A)$ to the antipodal map, also see the calculation in
the proof of \prettyref{prop:f} below.  With this fact the
introduction contains a complete, albeit sketchy, proof of
\prettyref{thm:lovasz}.
\end{rem}

\begin{defn}
We define $f_0\colon\HHom\to\Sphere^n$ as the composition of the
simplicial map 
\begin{align*}
\HHom&\to\esd(\bsd(\Delta^{n+1})),\\
\phi&\mapsto(\phi(v_0),\compl(\phi(v_{2r})))
\end{align*}
with the
triangulation~$t'$.
\end{defn}

\begin{lem}\label{lem:A0}
For any point $y\in\HHom$ we have $f_0(y)\ne f_1(y)$.
\end{lem}

\begin{proof}
Let $y\in\HHom$.  The point~$y$ 
is contained in the interior of a unique simplex
defined by a chain $\phi_0<\dots<\phi_s$.  The point $f_1(y)$ is in
the image under~$t$ of the simplex spanned by
$\phi_0(v_1)\subset\dots\subset\phi_s(v_1)$.  The point
$f_0(y)$ is in the image under $t'$ of the interior of the simplex
spanned by $(\phi_0(v_0),\compl\phi_0(v_{2r}))
\le\dots\le(\phi_s(v_0),\compl\phi_s(v_{2r}))$ and hence in the image
under~$t$ of the interior of the simplex spanned by
$\phi_0(v_0)\subset\dots\subset\phi_s(v_0)
\subset\compl\phi_s(v_{2r})\subset\dots\subset\compl\phi_0(v_{2r})$.
If $f_0(y)=f_1(y)$ then
$\phi_0(v_1)=\phi_0(v_0)$ contradicting
$\phi_0(v_1)\intersect\phi_0(v_0)=\emptyset$.
\end{proof}

\begin{defn}
We define $A_{r,n}\subset X_{r,n}$ by
\begin{align*}
A_{r,n}&\deq\Union_{i=0}^{2r-1} A_{r,n}^i,\\
A_{r,n}^i&\deq\set{x\in X_{r,n}\colon x_{i}=x_{i+1}},\qquad 0\le i<2r-1,\\
A_{r,n}^{2r-1}&\deq\set{x\in X_{r,n}\colon x_{2r-1}=-x_0}.
\end{align*}
Since $\tau A_{r,n}^i=A_{r,n}^{2r-1-i}$, the subset $A_{r,n}$ of $X_{r,n}$ is
invariant under the operation of~$\Z_2$.
\end{defn}

\begin{prop}\label{prop:f}
The above construction yields a $\Z_2$-equivariant map
\begin{equation*}
f\colon\HHom\xto\quad X_{n,r}\wo A_{n,r}.
\end{equation*}
\end{prop}

\begin{proof}
Let $y\in\HHom$.  We assume that $y$ has affine coordinates
$(\lambda_j)_j$ in the simplex $\simp{\phi_0,\dots,\phi_s}$ and write this fact
as $y=\sum_j \lambda_j\phi_j$.  For $i>0$ we have 
\begin{align*}
f_i(\tau y)
&=t\Biggl(\sum_j\lambda_j(\tau\phi_j)(v_{i})\Biggr)
=t\Biggl(\sum_j\lambda_j\phi_j(v_{2r-i})\Biggr)
\\&=f_{2r-i}(y).
\intertext{Also,} 
f_0(\tau y)
&=t'\Biggl(\sum_j\lambda_j((\tau\phi_j)(v_0),\compl(\tau\phi_j)(v_{2r}))\Biggr)
=t'\Biggl(\sum_j\lambda_j(\phi_j(v_{2r}),\compl\phi_j(v_0))\Biggr)
\\&=t\Biggl(\sum_j\lambda_j\left(
    \frac12\phi_j(v_{2r})+\frac12\compl\phi_j(v_0)\right)\Biggr)
\\&=-t\Biggl(\sum_j\lambda_j\left(
    \frac12\compl\phi_j(v_{2r})+\frac12\phi_j(v_0)\right)\Biggr)
=-t'\Biggl(\sum_j\lambda_j(\phi_j(v_0),\compl\phi_j(v_{2r}))\Biggr)
\\&=-f_0(y).
\end{align*}
Therefore $f(\tau y)=\tau f(y)$.

For $0<i<2r-1$ we have $f(y)\notin A_{r,n}^i$ by \prettyref{lem:Ai}.
From \prettyref{lem:A0} follows $f(y)\notin A_{r,n}^0$ and $f(y)=\tau
f(\tau y)\notin\tau A_{r,n}^0=A_{r,n}^{2r-1}$.  
Thus $f(y)\notin A_{r,n}$.
\end{proof}

\section{$\w_1^n(X_{r,n}\wo A_{r,n})=0$}
\label{sec:w}
We want to determine the class $\w_1^n(X_{r,n}\wo A_{r,n})\in
H^n((X_{r,n}\wo A_{r,n})/\Z_2; \Z_2)$.  Given any equivariant map
$X_{r,n}\wo A_{r,n}\to\Sphere^n$, this class is the image of the
generator $\omega\in H^n(\R P^n;\Z_2)$ under the induced map
$H^\ast(\R P^n;\Z_2)\to H^\ast((X_{r,n}\wo A_{r,n})/\Z_2; \Z_2)$.  The
equivariant map we will use is the projection
\begin{align*}
p_0\colon X_{r,n}=\prod_{i=0}^{2r-1}\Sphere^n&\to\Sphere^n\\
    x&\mapsto x_0.
\end{align*}
We note that the induced map $\bar p_0\colon X_{r,n}/\Z_2\to\R P^n$ is
a locally trivial fibration with fiber $(\Sphere^n)^{2r-1}$.

Since $A_{r,n}/\Z_2$ is a closed subset of the $2rn$-dimensional
manifold $X_{r,n}/\Z_2$, we have Poincaré duality at our disposal.  We
notationally omit the coefficient field~$\Z_2$ from now on and
consider the commutative diagram
\begin{equation*}\xymatrix@C-1em{
&&H^n(\R P^n)
\ar[dl]_{\bar p_0^\ast}
\ar[d]
\\
H^n(X_{r,n}/\Z_2,(X_{r,n}\wo A_{r,n})/\Z_2)
\ar[r]_-{k^\ast}
\ar[d]^\isom
&
H^n(X_{r,n}/\Z_2)\ar[r]_-{j^\ast}\ar[d]^\isom
&
H^n((X_{r,n}\wo A_{r,n})/\Z_2)
\\
H_{(2r-1)n}(A_{r,n}/\Z_2)
\ar[r]
&
H_{(2r-1)n}(X_{r,n}/\Z_2)
\\ 
H_{(2r-1)n}(A_{r,n})
\ar[u]
\ar[r]^-{l_\ast}
& H_{(2r-1)n}(X_{r,n})
\ar[u]^{\pi_\ast} }
\end{equation*} 
with $j$, $k$, and $l$ inclusions
und the isomorphisms given by Poincaré duality.

Now $\w_1^n(X_{r,n}\wo A_{r,n})=(\bar p_0\cmps j)^\ast(\omega)$, where
$\omega$ is the generator of $H^n(\R P)$.  The Poincaré dual of $\bar
p_0^\ast(\omega)$ is $\hbox{$\bar p_0$}_!(\ast)$ with $\ast$ being the
class of a point, i.e.\ the generator of $H_0(\R P^n)$, the Poincaré
dual of~$\omega$.  $\hbox{$\bar p_0$}_!(\ast)$ is the orientation
class (a bit of a misnomer when working with coefficients in~$\Z_2$)
of a fibre of the map $\bar p_0$, and therefore $\hbox{${{\bar
p}_0}$}_!(\ast)=\pi_\ast(c_0)$ with the following notation.
\begin{nota}
For $0\le i<2r$ we define $c_i\in H_{(2r-1)n}(X_{r,n})$ by
\begin{equation*}
c_i\deq \underbrace{[S^n]\times\dots\times[S^n]}_{\text{$i$ factors}}
  \times\ast\times
  \underbrace{[S^n]\times\cdots\times[S^n]}_{\text{$2r-i-1$ factors}}.
\end{equation*}
\end{nota}

\begin{lem}
For $0\le i<2r-1$ we have $c_i+c_{i+1}\in\im l_\ast$.
\end{lem}

\begin{proof}
Considering the map $\Delta\colon \Sphere^n\to\Sphere^n\times\Sphere^n$
given by
$\Delta(x)=(x,x)$ we have $\Delta_\ast(c)=\ast\times c+c\times\ast$
for $c\in H_n(\Sphere^n;\Z_2)$.  With
\begin{align*}
\Delta_i&\colon (\Sphere^n)^{2r-1}\to A_{r,n}^i\subset A_{r,n}\\
\Delta_i&\deq \underbrace{\id\times\dots\times\id}_{\text{$i$ factors}}
  \times\Delta\times
  \underbrace{\id\times\cdots\times\id}_{\text{$2r-i-1$ factors}}
\end{align*}
we therefore have
$c_i=l_\ast((\Delta_i)_\ast([\Sphere^n]\times\cdots\times[\Sphere^n]))$.
\end{proof}

\begin{lem}
We have $\pi_\ast(c_r)=0$.
\end{lem}

\begin{proof}
Singling out the sphere at position~$r$ of $X_{r,n}$ we obtain a
$\Z_2$-homeomorphism $X_{r,n}\homeo Y\times\Sphere^n$, where $Y$ is a
product of $2r-1$ spheres with the appropriate free $\Z_2$-action and
the operation on $\Sphere^n$ is trivial.  Under this homeomorphism,
$c_{r}$ corresponds to $[Y]\times\ast$, where $[Y]$ is the orientation
class of $Y$.  Since the quotient map $Y\to Y/\Z_2$ is of even degree,
the image of $[Y]\times\ast$ in
$H_{(2r-1)n}((Y\times\Sphere^n)/\Z_2)=H_{(2r-1)n}((Y/\Z_2)\times\Sphere^n)$
is zero.
\end{proof}

\begin{prop}\label{prop:wXA}
We have $w_1^n(X_{r,n}\wo A_{r,n})=0$.
\end{prop}

\begin{proof}
As noted before, $\w_1^n(X_{r,n}\wo A_{r,n})=j^\ast(\bar
p_0^\ast(\omega))$ with $\bar p_0^\ast(\omega)$ Poincaré dual to
$\pi_\ast(c_0)$.  By the preceding lemmas
\begin{equation*}
\pi_\ast(c_0)=\pi_\ast(c_0+c_r)=\pi_\ast
\left(\sum\nolimits_{i=0}^{r-1}(c_i+c_{i+1})\right)
  \in\im(\pi_\ast\cmps l_\ast).
\end{equation*}
Using the commutative diagram it follows that $\bar
p_0^\ast(\omega)\in\im k^\ast=\ker j^\ast$.
\end{proof}

\section{Result}

Our main result appears as \prettyref{thm:wv} in the introduction.
\begin{thm*}
Let $r\ge1$ and $n\ge0$.  Then $\w_1^n(\HHom)=0$.
\end{thm*}

\begin{proof}
This follows directly from \prettyref{prop:wXA}, since
$\w_1$ is characteristic: With the map $\bar f\colon\HHom/\Z_2\to
(X_{r,n}\wo A_{r,n})/\Z_2$ induced by the equivariant map~$f$ from
\prettyref{prop:f} we have
\[\w_1^n(\HHom)=\bar
f^\ast(\w_1^n(X_{r,n}\wo A_{r,n}))=\bar f^\ast(0)=0\]as claimed.
\end{proof}

\bibliographystyle{cas-ea}
\bibliography{math,topology,combi}

\end{document}